\documentstyle[12pt,psfig]{article}

\input{psfig}

\begin{document}

\title{\Large \bf
Surgeries on periodic links and homology of periodic 3-manifolds}

\author{ J\'ozef H.Przytycki and Maxim Sokolov }
\date{}

\maketitle

\begin{abstract}
We show that a closed orientable 3-manifold $M$ admits an action of
${\bf Z}_p$ with fixed point set $S^1$ iff $M$ can be obtained as the result of
surgery on a $p$-periodic framed link $L$ and ${\bf Z}_p$ acts freely
on the components of $L$.  We prove a similar theorem for free
${\bf Z}_p$-actions.  As an interesting application, we prove the
following, rather unexpected result:
for any $M$ as above and for any odd prime $p$, $H_1(M, {\bf Z}_p)\ne {\bf Z}_p$.
We also prove a similar criterion of 2-periodicity for rational
homology 3-spheres.
\end{abstract}

 \noindent {\large \bf 0. Introduction}

\medskip
In the early 1960's both Wallace [Wa] and Lickorish [Li] proved
that every closed, connected, orientable 3-manifold may be
obtained by surgery on a framed link in $S^3$. Thus, link diagrams
may be used to depict manifolds. Every manifold has infinitely
many different framed link descriptions. However, in
1970's Kirby [K] showed that two framed links determine the same
3-manifold iff they are related by a finite sequence of two
specific types of moves. This calculus of framed links, together
with the earlier results, gives a classification of 3-manifolds in
terms of equivalence classes of framed links. The framed link
representation of 3-manifolds has proven to be extremely useful.
For instance, most of the new 3-manifold invariants originating
from famous Witten's paper [Wi] are based on the framed link
approach. Therefore, it is always very useful to have some kind
of correspondence between certain classes of 3-manifolds and some
classes of framed links. One example of such correspondence is
the classical relationship between the lens spaces and the
chain-link diagrams (see, for instance, [Ro]). This result, in
particular, allowed L. Jeffrey to determine exact formulas for the
Witten-Reshetikhin-Turaev invariants of the lens spaces [Je].
Another example is the fact that a closed oriented 3-manifold is
an integral homology 3-sphere iff it can be obtained by surgery on
an algebraically split link with framing numbers $\pm 1$ (see
[Mu-1], [O-1]). This relationship plays a key role in many papers on
quantum and finite invariants of integral homology 3-spheres (see,
for instance, [O-1], [O-2], [Mu-1], [Mu-2]). In Section 1 we will establish an
analogous relationship between periodic 3-manifolds and periodic
links. Namely, we prove the following theorem:

\medskip
\noindent
{\sc Theorem 1.1} {\it Let $p$ be a prime integer and $M$ be a closed
oriented 3-manifold. There is an action of the cyclic group ${\bf Z}_p$
on $M$ with the fixed-point set equal to a circle if and
only if there exists a framed $p$-periodic link $L \subset S^3$ such that $M$ is the
result of surgery on $L$ and ${\bf Z}_p$ acts freely on the set of
components of $L$. }

\medskip
A special case of Theorem 1.1 when $M$ is a homology sphere was
proven in [Ka-Pr]. In the general form the theorem was proven for
the first time by the first author in his graduate course {\it
Topics in Algebra Situs}\footnote{The George Washington
University, February of 1999.}. A similar result is obtained for
manifolds with free ${\bf Z}_p$ actions:

\medskip
\noindent {\sc Theorem 1.2} {\it Let $p$ be a prime integer and
$M$ be a closed oriented 3-manifold. There is a free action of the
cyclic group ${\bf Z}_p$ on $M$ iff there exists a framed $p$-periodic
link $L \subset S^3$ admitting a free action of ${\bf Z}_p$ on the
set of its components such that $M$ is the result of surgery on
$L'= L\cup \gamma$, where $\gamma$ is the axis of the action with
framing co-prime to $p$. }

\medskip
In Section 2 we give an interesting application of Theorem 1.1.
Namely, we prove the following result:

\medskip
\noindent
{\sc Theorem 2.1}
{\it If a closed orientable 3-manifold $M$ admits an action of a
cyclic group ${\bf Z}_p$ where $p$ is an odd prime integer and the
fixed point set of the action is $S^1$ then $H_1(M; {\bf Z}_p) \ne
{\bf Z}_p$.}

\medskip
Note that this theorem provides a non-trivial criterion for
3-manifolds admitting the described action. The simplest examples
of 3-manifolds with $H_1(M; {\bf Z}_p) = {\bf Z}_p$ are lens
spaces $L_{pn, q}$ (or more generally, $(pn/q)$ Dehn surgeries on
knots in $S^3$).

Theorem 2.1 was first announced as a conjecture\footnote{ The
conjecture was obtained as a result of extensive computations
performed with a program written in {\it Mathematica}.}
 and
partially proven\footnote{In the case when the orbit space of the
action can be obtained from $S^3$ by an integer surgery on a knot. }
in April of 1999 [So] (It is interesting to mention
that the conjecture was influenced by the study of
Murakami-Ohtsuki-Okada invariants on periodic 3-manifolds\footnote{
In turn, our interest in Murakami-Ohtsuki-Okada
invariants was sparked by their relation with the second skein module [P-4].}
, but the
equation $MOO_p(M)= \pm G_p^{rkH_1(M; {\bf Z}_p)}$ eventually led
to the more ``classical" algebraic topology.  Here $p$ is an odd prime
integer, $MOO_p$ is the Murakami-Ohtsuki-Okada invariant
parameterized by $q=e^{2 \pi i / p}$, and $G_p = \sum _{j\in {\bf Z}_p}
q^{j^2}$
). Recently (November,
1999), Adam Sikora announced a proof of the theorem
[Si]. In fact, using some classical but involved algebraic
topology, he obtained more general results implying our theorem.
The surgery presentation of periodic 3-manifolds
developed in Section 1 allowed us to find an elementary proof of
Theorem 2.1, presented in Section 2.

Theorem 2.1 is not true for $p=2$, see Remark 2.12.  An interesting
criterion for 2-periodic rational homology spheres is provided by
the following theorem.

\medskip
\noindent
{\sc Theorem 2.2.} {\it Let $M$ be a rational homology 3-sphere such
that the group $H_1(M; {\bf Z})$ does not have elements of order
16.  If $M$ admits an orientation preserving action of ${\bf Z}_2$
with the fixed point set being a circle then the canonical
decomposition of the group $H_1(M; {\bf Z})$ has even number of
terms ${\bf Z}_2$ and even number of terms ${\bf Z}_4$, and
arbitrary number of terms ${\bf Z}_8$. }

\medskip
The second author thanks Yongwu Rong and Adam Sikora for useful
conversations.  When a preliminary version of the paper was ready,
we received an e-mail from James Davis saying that he and his
student Karl Bloch also found a proof for Theorem 2.1.

\ \\

\bigskip
 \noindent {\large \bf 1. Periodic 3-manifolds are surgeries
on periodic links}

\medskip
We show in this section that $p$-periodic closed oriented
3-manifolds can be presented as results of integer surgeries on
$p$-periodic links. We show also an analogous result for manifolds
with free action of ${\bf Z}_p$.

Before we prove theorems 1.1 and 1.2, we need to establish some basic
terminology and preliminary lemmas. \ \\ \

\noindent
{\bf 1.1 Periodic Links.}

\medskip
\noindent
{\sc Definition.} By a {\it framed knot} $K$ we mean a ring
$S^1\times [ 0, \varepsilon ]$ embedded in $S^3$. By the {\it
framing} of $K$ we mean an integer defined as follows.  Let
$V_{\varepsilon}$ be the $\varepsilon$-neighborhood of
$K_0=S^1\times \{ 0\}$, then $K_{\varepsilon}= S^1 \times
\{\varepsilon\}$ is a projection of $K_0$ onto $\partial
V_{\varepsilon}$.  Let $P$ be the projection of $K_0$ onto $\partial
V_{\varepsilon}$which is homologically trivial in $S^3 -
V_{\varepsilon}$.  The framing $f$ is defined as the algebraic
number of intersections of $K_{\varepsilon}$ and $P$.  A {\it
framed link} is a collection of non-intersecting framed knots.  We
will adopt the usual ``blackboard" convention for framed link
diagrams.

\medskip

\noindent {\sc Definition.} A (framed) link $L$ in $S^3$ is called
{\it $p$-periodic} if there is a ${\bf Z}_p$-action on $S^3$, with a circle as
a fixed point set, which maps $L$ onto itself, and such that $L$
is disjoint from the fixed point set. Furthermore, if $L$ is an
oriented link, one assumes that each generator of ${\bf Z}_p$ preserves
the orientation of $L$ or changes it to the opposite one.

\medskip

By the positive solution of Smith Conjecture ([M-B],
[Th]) we know that the fixed point set of the action of ${\bf Z}_p$
is an unknotted circle and the action is conjugate to an
orthogonal action on $S^3$. In other words, if we identify $S^3$
with ${\bf R}^3\cup\infty$, then the fixed point set can be assumed to
be equal to the ``vertical'' axis $z=0$ together with $\infty$,
and a generator $\varphi$ of ${\bf Z}_p$ can be assumed to be the
rotation $\varphi(z,t) = (e^{{2\pi i}/{p}}\cdot z,t)$, where the
coordinates on ${\bf R}^3$ come from the product of the complex plane
and the real line ${\bf C}\times {\bf R}$. Thus, any (framed) $p$-periodic
link $L^p$ may be represented by a $\varphi$-invariant diagram,
{\it $p$-periodic diagram} (with framing parallel to the
projection plane), see Fig. 1.

\ \\
\centerline{\psfig{figure=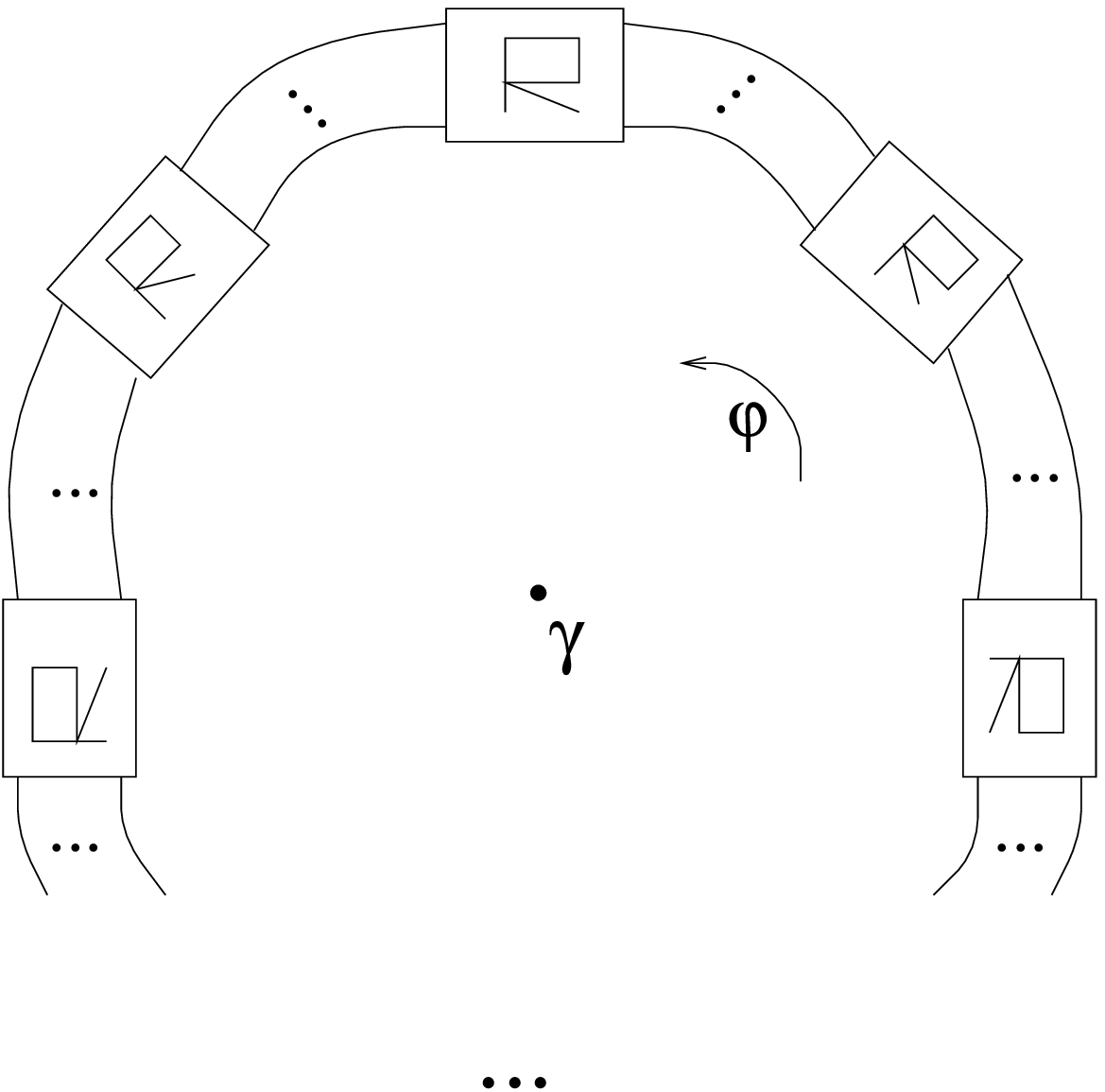,height=3.2cm}}
\begin{center}
Fig. 1
\end{center}
\

\noindent
By the {\it underlying link} for $L^p$ we will mean the
orbit space of the action, that is $L_* = L^p/ {\bf Z}_p$.

\medskip
\noindent
{\sc Lemma 1.3} {\it Let $p$ be a prime integer and
$L^p\subset S^3$ be a (framed) $p$-periodic
link. The following three conditions are equivalent:

\noindent
1) ${\bf Z}_p$ acts freely on the set of components of $L^p$;

\noindent
2) The linking number of each component of the underlying link
$L_*$ with the axis of rotation is congruent to zero modulo $p$;

\noindent
3) The number of components of $L^p$ is $p$ times greater
than the number of components of $L_*$.
}

\noindent
{\sc Proof.}
The equivalence of the conditions 1) and 3) is obvious.  To prove
that 2) is equivalent to 3), consider the covering projection
$L^p\to L_*$.  Let $l$ be a component of $L_*$.  The preimage $\rho ^{-1}$ of
the closed path $\lambda$ which traverses $l$ exactly once (i.e., $l$ with
a base point) consists of
$p$ paths $ \lambda _1, \dots , \lambda _p$ in $L^p$. The condition
 $\mbox{lk}(l, \gamma) =0 \pmod{p}$ is equivalent to the condition that
 each of $\lambda_i$ is closed. Thus, $l$ lifts to $p$
components in $L^p$ iff $\mbox{lk}(l, \gamma) =0 \pmod{p}$.
\hfill $\Box$

\medskip
\noindent
 {\sc Definition.}  A $p$-periodic link $L^p$ that satisfies any of the
conditions from Lemma 1.3 will be called {\it strongly
$p$-periodic}.

\ \\ \
{\bf 1.2 Periodic Manifolds.}

\noindent {\sc Definition.} A $3$-manifold $M$ is called {\it
$p$-periodic} if it admits an orientation preserving
action of the cyclic group ${\bf Z}_p$
with a circle as a fixed point set, and the action is free outside
the circle.

We can show immediately the easy part of Theorem 1.1. Indeed,
consider a strongly $p$-periodic framed link $L^p$, and let $M$ be
the 3-manifold obtained by surgery on $L^p$.  By definition of a
framed $p$-periodic link, there is a ${\bf Z}_p$ action on $S^3$,
and on $S^3-L^p$, with a circle $\gamma$ as a fixed point set.
This action induces a ${\bf Z}_p$ action on $M$. Moreover, since
the action of ${\bf Z}_p$ is free on the set of components of
$L^p$, there are no other fixed points of the action of ${\bf
Z}_p$ on $M$ but the circle $\gamma$.

To show the difficult part of Theorem 1.1 we first
fix some notation.
 Suppose that ${\bf Z}_p$ acts on $M$ with the
fixed-point set equal to a circle $\gamma$. Denote the quotient by
$M_*=M/{\bf Z}_p$, the projection map by $h \colon M\to M_*$ and
$\gamma_*=h(\gamma)$.

\medskip
\noindent {\sc Lemma 1.4} {\it The map $h_* \colon H_1(M)\rightarrow
H_1(M_*)$ is an epimorphism. }

\noindent
 {\sc Proof.}
Let $x_0\in\gamma$. Since $x_0$ is a fixed point of the action,
any loop based at $h(x_0)$ lifts to a loop based at $x_0$. Thus
$h_{\# }: \pi_1(M,x_0)\rightarrow \pi_1(M_*,h(x_0))$ is an
epimorphism, and since $H_1$ is an abelianization of $\pi_1$, the
map $h_*$ is also an epimorphism. \hfill $\Box$

\medskip
Notice that the proof works for any finite group action on a manifold
with a non-empty fixed point set.

Let us recall the Lefschetz' duality theorem which we will use in
our proof of Theorem 1.1. First some terminology: a compact
connected $n$-dimensional manifold $M$ is called {\it $R$-oriented} for
a commutative ring with identity $R$, if $H_n(M,\partial M;R)=R$.
In particular, any manifold is ${\bf Z}_2$-oriented, and an oriented
manifold is $R$-oriented for any ring $R$. For a reference, see
[Sp].

\medskip
\noindent {\sc Theorem 1.5} (Lefschetz) {\it Let $M$ be a compact
$n$-dimensional, $R$-oriented manifold. Then there is an
isomorphism $\tau : H^q(M;R) \to H_{n-q}(M,\partial M;R)$.
Furthermore if $R$ is a PID (principal ideal domain) and
$H_{q-1}(M,R)$ is free then $H^q(M;R)= Hom(H_q(M;R), R)$ and for
$\alpha \in H^q(M;R)$ and $c\in H_q(M;R)$ one has: $\alpha (c) =
alg(c,\tau(\alpha))$, where $alg(c,\tau(c))\in R$ is the algebraic
intersection number of $c$ and $\tau(\alpha)$ in $M$ ( $alg \colon
H_q(M;R) \times H_{n-q}(M,\partial M;\\ R) \to R$).}

\medskip
We use Lefschetz' Theorem to show that the covering $h: M \to M_*$
is yielded by a 2-chain whose boundary is a multiple of
$\gamma_*$. Because we work with $q=1$, then $H_{q-1}(M,R)$
is free and we can use the intersection number interpretation of
the Lefschetz' Theorem.

\medskip
\noindent
{\sc Lemma 1.6} {\it
Let $M$ be a closed orientable $p$-periodic 3-manifold.
With the notation as before, one has:
\begin{enumerate}
\item[(1)]
$\gamma_* \equiv 0$ in $H_1(M_*,{\bf Z}_p)$.
\item[(2)]
There is a 2-chain $C\in C_2 (M_*, {\bf Z}_p)$ such that $\partial
C \equiv m\gamma_* \ mod \ p$ and the covering $h: (M-\gamma) \to
(M_* - \gamma_*)$ is yielded by the map $\phi_C : H_1(M_* -
\gamma_*) \to {\bf Z}_p$ where $\phi_C(K)$ is the intersection number of
$K$ with $C$ (i.e. for a 1-cycle $K \in M_* - \gamma_*$,
$\phi_C(K)= alg(K,C)$ where $alg(K,C)$ is the intersection number
of $K$ with $C$, well defined $mod \ p$)\footnote{
If $H_2(M_*, {\bf Z})=0$, then $alg(K,C) = lk(K,m\gamma_*)$,
but generally $alg(K,C)$ depends on the choice of $C$.}.
In particular
$\phi_C(\mu_*)= m$, where $\mu_*$ is a meridian of $\gamma_*$.
\end{enumerate}
}

\noindent
 {\sc Proof.}
  To work with Lefschetz' Theorem we have to consider
compact manifolds. Thus, instead of $M_* - \gamma _*$ we consider
a homotopically equivalent compact manifold $\hat M_* = M_*
-int(V_{\gamma_*})$, where $V_{\gamma_*}$ is a regular
neighborhood of $\gamma_*$ in $M_*$. Similarly, let $V_{\gamma} =
h^{-1}(V_{\gamma_*})$ be a $Z_p$-invariant regular neighborhood of
$\gamma$ in $M$. Let also $\hat M = M- int(V_{\gamma})$.
 Since $\hat h \colon \hat M \to \hat M_*$ is a regular covering,
 it is characterized by an epimorphism $\pi _1(\hat M_*) \to
  \pi _1(\hat M_*)/ \pi _1(\hat M) = {\bf Z}_p$ (up to an
  automorphism of ${\bf Z}_p$).  Thus, since ${\bf Z}_p$ is
  abelian, $\hat h$ is defined by an epimorphism  $\phi : H_1(\hat M_*) \to
{\bf Z}_p$, where $\phi$ is unique up to an automorphism of ${\bf Z}_p$.
 Let $\hat
C$ be a 2-cycle representing the element of $H_2(\hat M_*,\partial
\hat M_*; {\bf Z}_p)$ dual to the epimorphism $\phi$, that is, such that
$alg(K,\hat C) \equiv \phi(K)\pmod{p}$ for any $K\in H_1(\hat M_*)$.
We can assume that $\hat C \cap
\partial\hat M_*$ is a collection of simple noncontractible curves
in the torus $\partial \hat M_*$. Finally let $C$ be a 2-chain
obtained from $\hat C$ by adding to $\hat C$ annuli connecting the
components of $\hat C \cap \partial\hat M_*$ with $\gamma_*$ in
$V_{\gamma_*}$. Thus, $C$ of part (2) is constructed.

Let $\mu_*$ be a meridian of $\gamma_*$ (or more precisely, of
$V_{\gamma_*}$). $\phi (\mu_*) \neq 0\pmod{p}$ because $\gamma_*$ is a
branching set of the covering, so the preimage of $\mu_*$ (under
$h$) is a connected curve $\mu$, (a meridian of $\gamma$ in $M$),
by the definition of a branched covering. Thus, there is $0<m<p$ such
that  $\phi(\mu_*) = m$. We can conclude also that $m\gamma_* \equiv \partial C$
mod $p$, thus $\gamma_* \equiv 0$ in $H_1(M_*;{\bf Z}_p)$.
 \hfill $\Box$

\medskip
\noindent
{\sc Proof of Theorem 1.1.}
By the classical result of Wallace (1960) and Lickorish (1962),
every closed oriented 3-manifold is a result of a surgery on a
framed link in $S^3$. In particular $M_*$ can be represented as a
result of surgery on some framed link $L_\#$ in $S^3$. Inversely
$S^3$ can be obtained as a result of surgery on some framed link
$\hat L_\#$ in $M_*$. We can assume that $\hat L_\#$ satisfies the
following conditions (possibly after deforming $\hat L_\#$ by
ambient isotopy):
\begin{enumerate}

\item[(1)] $\gamma_*\cap \hat L_\#=\emptyset$;
\item [(2)] ${alg}(\hat L_\#^i, C)\equiv 0\ mod\ p$,
for any component $\hat L_\#^i$ of $\hat L_\#$.
\end{enumerate}
$\hat L_\#$ satisfying the conditions 1-2 can be obtained as
follows:\\ Let $L_\#\subset S^3$ be a framed link in $S^3$ such
that $M_*$ is a result of surgery on $L_\# $. Let $\hat L_\#$ denote
the co-core of the surgery \footnote{In our notation, {\it core}
of the surgery is the framed surgery link, that is the framed
link, regular neighborhood of which is removed in the ``drilling''
part of the surgery, with framing yielded by the meridian of the
attached (``filling" part of the surgery) solid torus. The {\it
co-core} of the surgery is the core of the ``filling" solid torus,
with its framing yielded by the meridian of the removed solid
torus. The surgery on the co-core link brings back the initial
manifold.}. In particular, $\hat L_\#$ is a framed link in $M_*$
such that $S^3$ is a result of surgery on $\hat L_\#$. By a general
position argument, we can make $\gamma_*$ and $\hat L_\#$ disjoint,
but in order to get condition (2) we should do so in a
controllable manner.

Let $C$ be the 2-chain from Lemma 1.6. Let $\hat L^i_\#$ be any
component of $\hat L_\#$. If we change a crossing between
$\hat L^i_\#$ and $\gamma_*$ then the algebraic crossing number,
$alg(\hat L^i_\#, C)$ changes by $\pm m$ mod $p$. Thus, by a
series of crossing changes we can get $alg(\hat L^i_\#,
C)\equiv 0 \ mod \ p$ for any component of $\hat L_\#$,
providing condition (2).

This implies that we can easily modify
$C$ (outside $\gamma_*$) so that $C \cap \hat L_\#= \emptyset$.
Therefore, $C$ survives the surgery (as well as $\gamma_*$), and in
$S^3$ it has ${\bf Z}_p$ boundary $m\gamma_*$ and it is disjoint from
$L_\#$ (link in $S^3$ being the co-core of the surgery on $\hat
L_\#$ in $M_*$). Thus, ${\rm lk}(L_\#^i, \gamma_*)\equiv 0\ mod\ p$
for any component $L_\#^i$ of $L_\#$.

Now we are ready to unknot $\gamma_*$ using Kirby calculus
([K], [F-R]). Choose some orientation on $\gamma_*$. We
can add unlinked components with framing $\pm 1$ to $L_\#$ around
each crossing of $\gamma_*$, making sure that arrows on $\gamma_*$
run opposite ways (i.e. the linking number of $\gamma_*$ with the
new component is zero, see Fig. 2). Use the K-move to
change the appropriate crossings and thus to unknot $\gamma_*$.
Thus we trivialized $\gamma_*$ without compromising conditions (1)
and (2). Denote the framed link obtained from the initial link $L_\#$
after the described isotopy and adding the new components by $L_*$.
 Notice, that each new component that we introduce during
the above procedure has linking number $0$ with any other component of
$L_*$.

\ \\
\centerline{\psfig{figure=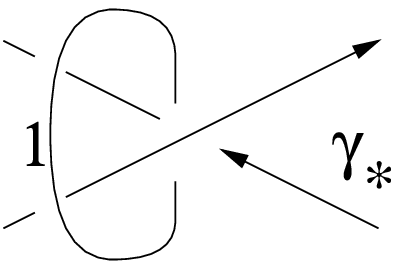}}
\begin{center}
Fig. 2
\end{center}

To complete the proof of the theorem, consider the $p$-fold cyclic
branched covering of $S^3$ by $S^3$ with branching set $\gamma_*$. Let $L$
denote the preimage of $L_*$. Notice that $L$ is strongly
$p$-periodic. We claim that the result of performing surgery on
$L$ is ${\bf Z}_p$-homeomorphic to $M$.

The preimage of each component of $L_*$ consists of $p$ components
permuted by a ${\bf Z}_p$ action, by Lemma 1.4. Therefore, ${\bf Z}_p$ acts
on the result of the surgery on $S^3$ along $L$, $\widetilde M=(S^3,L)$, with a
branched set $\gamma$ and quotient $M_*=(S^3,L_*)$.

The group $H_1(M_*-\gamma_* - \hat L_\#)=H_1(S^3-\gamma_* - L_*)$
is generated by $\mu_*$ (a
meridian of $\gamma_*$) and meridians of the components of $L_*$, say
$b_1,\dots ,b_k$.
Of course, $H_1(M_*-\gamma_*)$ is also generated by $\mu_*, b_1,\dots
,b_k$.  By Lemma 1.6, the covering $\rho \colon
(M-\gamma) \to (M_*-\gamma_*)$
is
characterized by the map  $\phi_C \colon H_1(M_*-\gamma_*)\to {\bf
Z}_p$ (up to an automorphism of ${\bf Z}_p$), where $\phi_C (K)$
is the intersection number  of $K\in H_1(M_*-\gamma_*) $ with
$C$ modulo $p$.
Similarly, the covering
 $\widetilde \rho \colon
(\widetilde M-\gamma) \to (M_*-\gamma_*)$ is characterized by a
map $\phi_2 \colon H_1(M_*-\gamma_*)\to {\bf Z}_p$.  By our
construction, $\phi_C(\mu _*)=m$ and $\phi _C(b_i)=0$ for every
$1\le i \le k$.  We need to show that $\phi_2(\mu_*)=m'$ for some
$m'$ coprime to p, and $\phi_2(b_i)=0$ for every $1\le i\le k$.
This follows from the fact that $\widetilde \rho ^{-1} (b_i)$
consists of $p$ loops and $\widetilde \rho ^{-1} (\mu_*)$ is a
single loop.  Thus, $\phi_2$ and $\phi_C$ are equivalent up to the
automorphism of ${\bf Z}_p$ sending $m'$ to $m$.  Therefore, the
manifolds $(\widetilde M -\gamma)$ and $(M-\gamma)$ are ${\bf
Z}_p$-homeomorphic, with a homeomorphism given by
$g\colon (\widetilde M -\gamma)\to (M-\gamma)$ such that $\widetilde
\rho = \rho \circ g$.   Notice, that $\widetilde
M$ can be obtained from $(\widetilde M -\gamma)$ by
attaching a 2-handle along $\widetilde \rho^{-1}(\mu_*)$ and then
a 3-handle, and
$M$ can be obtained from $(M-\gamma)$ by
attaching a 2-handle along $\rho^{-1}(\mu_*)$ and then
a 3-handle.  Since $g(\widetilde
\rho^{-1}(\mu_*))= \rho^{-1}(\mu_*)$, the homeomorphism $g$ can be
extended to a ${\bf Z}_p$-homeomorphism $\hat g\colon \widetilde
M\to M$.
Our
proof of the Theorem 1.1 is completed. \hfill $\Box$

\medskip
\noindent
{\sc Remark 1.7.} Notice, that if in the proof of Theorem 1.1 we
assumed that the link $L_\#$ was algebraically split then the link
$L_*$ would be algebraically split as well.  This remark will be
important later in the proofs of Theorems 2.1 and 2.2.

\medskip
As a corollary we obtain a proof of Theorem 1.2.

\medskip
\noindent {\sc Proof of Theorem 1.2.} Consider any ${\bf Z}_p$
equivariant knot, say $\hat \gamma$ in $M$. Let $V_{\hat \gamma}$
be a ${\bf Z}_p$ equivariant regular neighborhood of $\hat \gamma$ in
$M$ and $\gamma '$ a curve on $\partial V_{\hat \gamma}$ which is
also ${\bf Z}_p$ equivariant.  Notice, that $\gamma '$ intersects a
meridian of $V_{\hat \gamma}$ exactly once.
Now let $M'$ be a manifold obtained from
$M$ by a surgery on  $\hat \gamma$ with the framing defined by $\gamma
'$. Let $\gamma \in M'$ be the co-core of the surgery. The ${\bf Z}_p$
action on $M$ yields the action on $M'$ and our choice of framing
guarantees that $\gamma $ is the (only) fixed point set of the
action. Thus, we can apply to $M'$ the previous theorem.
This proves that $M$ can be obtained by an integer surgery on
$L\cup \gamma$, where $L$ is a strongly $p$-periodic link.
Furthermore, the framing of $\gamma$ must be coprime to $p$, to
insure that ${\bf Z}_p$ acts on the resulting manifold with no
fixed points.
\hfill
$\Box$

\medskip
\noindent
{\sc Remark 1.8.}
We plan to extend Theorem 1.1 to any ${\bf Z}_p$ orientation preserving
action on a closed 3-manifold $M$, and to ${\bf Z}_{p^k}$ actions.

\ \\

\noindent
{\large \bf 2. Homology of periodic 3-manifolds.}

\medskip
The main goal of this section is to give elementary proofs of
Theorems 2.1 and 2.2 using the surgery presentation of $p$-periodic
3-manifolds developed in Section 1.

\bigskip
\noindent {\bf 2.1 Linking matrices of framed strongly
$p$-periodic links and of algebraically split links.}

\medskip
Let $L$ be a framed oriented link of $n$ components $l_1, \dots ,
l_n$. The {\it Linking matrix} of $L$ is the matrix
$(a_{ij})_{n\times n}$ defined by

$$
 a_{ij} =
 \cases{
 \mbox{lk}(l_i, l_j) & if $i\ne j$\cr
\mbox{framing of }l_i & if $i=j$
}
$$

Let $L^p$ be a framed strongly $p$-periodic link and $L_*$ be the
corresponding underlying link. Fix an orientation of $L_*$ and
denote the components of $L_*$ by $l_1,\dots , l_n$. Consider a
$p$-periodic diagram of $L^p$. Denote the $p$ copies of the
tangle $R$ from the diagram (see Fig. 1) by $R_1, \dots , R_p$ in
the clockwise order. Lift the orientation of $L_*$ to $L^p$. By
Lemma 1.3, each component $l_i$ of $L_*$ has $p$ covering
preimages in $L^p$. Denote them in a clockwise order by $l_{i1},
\dots , l_{ip}$. By the clockwise order here we mean such an
order that if we choose any point $x\in l_{ij} \cap R_j$ then the
corresponding point in $R_{j+1}$ will belong to $l_{i j+1}$
(subscripts are treated modulo $p$).

Now consider the following natural order for the components of
$L^p$: $$
l_{11},\dots , l_{1p}, l_{21}, \dots ,l_{2p},\dots ,
l_{n1},\dots , l_{np}.
$$
It is not hard to see that with regard
to this order, the linking matrix $A_p$ for $L^p$ is of the
following form $$
A_p = \pmatrix{ A_{11} & B_{12} & \cdots & B_{1n}\cr
B_{21} & A_{22} & \cdots & B_{2n}\cr
\vdots & \vdots & \ddots & \vdots \cr
B_{n1} & B_{n2} & \cdots & A_{nn}\cr},
$$ where all the blocks are $p\times p$, $A_{ii}$ is the linking
matrix for the sublink consisting of $l_{i1}, \dots , l_{ip}$, and
$B_{ij}$ is the matrix with elements $b_{ks}^{ij} = \mbox{lk}
(l_{ik}, l_{js})$.

Recall, that a matrix $(a_{ij})_{k\times k}$ is called {\it
circulant} if $a_{ij} = a_{i+1 j+1}$, $i,j=1,\dots ,k$ (subscripts
mod $k$).

\medskip
\noindent
{\sc Proposition 2.3.} {\it Every block in $A_p$ is a
circulant matrix.}

\noindent {\sc Proof.} Consider $B_{ij} = (b_{ks}^{ij})_{p\times
p}$. Then $b_{ks}^{ij} = \mbox{lk} (l_{ik}, l_{js})$ and $b_{k+1
s+1}^{ij} = \mbox{lk} (l_{i k+1}, l_{j s+1})$. If one rotates the
$p$-periodic diagram of $L^p$ around the center in the clockwise
direction by $2\pi /p$ then the pair $(l_{ik}, l_{js})$ will go
into $(l_{i k+1}, l_{j s+1})$, taking subscripts modulo $p$. Thus,
$\mbox{lk}(l_{ik}, l_{js}) =\mbox{lk} (l_{i k+1}, l_{j s+1})$. If
we notice that the framing numbers of $l_{i1}, \dots , l_{ip}$ are
all the same, then the above argument shows that $A_{ii}$ is also
circulant for any $i=1,\dots , n$. \hfill $ \Box$

\medskip
\noindent {\sc Definition.} We will call a (framed) link $L$ {\it
algebraically split} if the linking number between any two
components of $L$ is zero. A strongly $p$-periodic (framed) link
$L^p$ will be called {\it orbitally separated} if the underlying
link $L_*$ is algebraically split.

\medskip
\noindent {\sc Remark 2.4.} It is not hard to see that a strongly
$p$-periodic link $L^p$ is orbitally separated iff any two
components of $L^p$ that cover different components of $L_*$ have
the linking number equal to zero.

\medskip
\noindent {\sc Corollary 2.5.} {\it It follows from Proposition
2.3 and Remark 2.4 that $L^p$ is an orbitally separated link iff
all the non-diagonal blocks $B_{ij}$ in $A_p$ are zero matrices.}
\hfill $\Box$

\bigskip
\noindent
{\bf 2.2 Nullity of symmetric circulant matrices over ${\bf Z}_p$.}
\medskip

Circulant matrices are very well studied and a lot is known about
them (see, for instance, [D]). But, apparently, not much is known
about circulant matrices over finite fields (or rings). The
following two results provide the key tool for our proof of
Theorem 2.1, but they also appear to be interesting from a purely
matrix theoretical point of view.

\medskip
\noindent
{\sc Lemma 2.6.} {\it If
$$
A= \pmatrix{ a_1 & a_2 & a_3 &
\cdots & a_n\cr
a_n & a_1 & a_2 & \cdots & a_{n-1}\cr
\vdots & \vdots &\vdots & \ddots & \vdots \cr
a_2 & a_3 & a_4& \cdots & a_1\cr}
$$
is a circulant matrix with integer elements then
$$
 \mbox{det }A =
 \cases{
a_1^n + a_2^n + \dots + a_n^n \pmod{n}, & if n is odd;\cr
a_1^n - a_2^n + \dots - a_n^n \pmod{n}, & if n is even.
}
$$
 }

\noindent {\sc Proof.} The determinant of $A$ is a sum of $n!$
terms. The terms of the form $a_i^n$, $i=1,\dots ,n$, will be
called {\it diagonal}. Note that any term different from diagonal
appears in the sum exactly $n$ times: $$ a_{i_1} a_{i_2} a_{i_3}
\dots a_{i_n}, $$ $$ a_{i_2} a_{i_3} \dots a_{i_n} a_{i_1}, $$ $$
\vdots $$ $$ a_{i_n} a_{i_1} a_{i_2} \dots a_{i_{n-1}}. $$ Note
that the sign for all such terms is the same. To prove this we
need to show that the permutations $$\sigma =\pmatrix{ 1 & 2 & 3 &
\cdots & n\cr i_1 & i_2 & i_3 & \cdots & i_n\cr} \mbox{ and }
\sigma '=\pmatrix{ 1 & 2 & 3 & \cdots & n\cr i_n+1 & i_1+1 & i_2+1
& \cdots & i_{n-1}+1\cr} $$ have the same parity (everything is
modulo $n$). Obviously,
\begin{eqnarray*}
\pmatrix{ 1 & 2 & 3 & \cdots & n\cr
i_n+1 & i_1+1 & i_2+1 & \cdots & i_{n-1}+1\cr} = \\
\pmatrix{ 2 & 3 & \cdots & n & 1\cr
i_1+1 & i_2+1 & \cdots & i_{n-1}+1 & i_n+1\cr}.
\end{eqnarray*}
The row $\pmatrix{ 2 &3 &4 &\cdots &n &1}$ has $n-1$ inversions.
The numbers of inversions in $\pmatrix{ i_1 &i_2 &\cdots &i_n}$
and $\pmatrix{ i_1+1 & i_2+1 & \cdots & i_n+1}$ also differ by
$n-1$. Thus the parities of $\sigma$ and $\sigma '$ are the
same. Therefore, the total sum of all non-diagonal terms in
$\mbox{det } A$ is 0 mod $n$. The result follows. \hfill $\Box$

\medskip
Denote the nullity of $A$ over ${\bf Z}_n$ by
 $\mbox{\it null}_n A$.

\medskip
\noindent
{\sc Lemma 2.7.} {\it Let $p$ be an odd prime integer and $A$ be
 a $p\times p$ symmetric
circulant matrix over ${\bf Z}_p$, then
$\mbox{null}_p A\ne 1$.}

\noindent {\sc Proof.} Assume $\mbox{\it null}_p A>0$, i. e.
$a_1^p + 2 a_2^p + \dots + 2 a_{ {{p+1}\over 2}}^p =
a_1 + 2 a_2 + \dots + 2 a_{ {{p+1}\over 2}} = 0 \pmod{p}$,
by Lemma 2.6 (the second equality follows from Fermat's theorem).
After adding all rows to the last one and all
columns to the last column we get
 $$ \mbox{\it det } A = \mbox{\it
det} \pmatrix{ a_1 & a_2 & a_3 & \cdots & a_3 & 0\cr a_2 & a_1 &
a_2 & \cdots & a_4 & 0\cr \vdots & \vdots &\vdots & \ddots &
\vdots \cr a_3 & a_4 & a_5& \cdots & a_1 & 0 \cr 0 & 0 & 0 &
\cdots & 0 & 0\cr}. $$
 Denote the $i$th column of the above matrix
by $C_i$. Then the linear combination
 $$ (p-1) C_1 + (p-2) C_2 +
\dots +2 C_{p-2} + C_{p-1}
 $$
  is 0 modulo $p$. Indeed, the $i$th
row of the linear combination is
 $$ (p-1) a_i + (p-2) a_{i+1} +
\dots + 2 a_{i-3} + a_{i-2},
 $$
  all the coefficients and
subscripts are modulo $p$. It is not hard to see that after the
substitution $a_1=-2 a_2 - 2 a_3 - \dots - 2 a_{ {{p+1}\over
2}}$, the coefficient for $a_k$ ($ k\ne1$) in the above sum is
 $$
-2(p-i) + (p-(i+k)) + (p-(i-k)) =0 .
$$
 Thus, if $\mbox{\it det }
A=0$ over ${\bf Z}_p$ then $\mbox{\it null}_p A\ge 2$. \hfill
$\Box$

\medskip
\noindent
{\sc Remark 2.8.} Lemma 2.7 is not true if $p=2$.
For instance, nullity of $\pmatrix{ 1 & 1\cr
1 & 1\cr}$
is 1 over ${\bf Z}_2$.

\bigskip
\noindent
{\bf 2.3. Proof of Theorem 2.1.}

\medskip
In this section we will prove Theorem 2.1, the main theorem of
Section~2. Let $M$ be a closed oriented 3-manifold obtained by a
Dehn surgery on a framed oriented link $L$, and let $A$ be the
linking matrix of $L$. The following fact is well-known.

\medskip
\noindent
{\sc Lemma 2.9.} {\it
$\mbox{null}_p A= \mbox{rank }H_1(M; {\bf Z}_p)$.} \hfill $\Box$

\medskip

Now we are ready to prove an important special case of Theorem
2.1.

\medskip
\noindent {\sc Proposition 2.10.} {\it Let $p$ be an odd prime integer.
If a closed orientable 3-manifold $M$ can be obtained from $S^3$
by Dehn surgery on an orbitally separated framed link $L^p$ then
$H_1(M; {\bf Z}_p)\ne {\bf Z}_p$.}

\medskip
\noindent {\sc Proof.} Assume that $M$ can be obtained by Dehn
surgery on an orbitally separated framed link $L^p$. Let $A_p$ be
the linking matrix of $L^p$, as constructed in Section 2.1. By
Corollary 2.5, $A_p$ is a block diagonal matrix. Therefore,
$\mbox{\it null}_p A_p$ is equal to the sum $\mbox{\it null}_p
A_{11} + \dots + \mbox{\it null}_p A_{nn}$. By Proposition 2.3,
each $A_{ii}$ is a circulant matrix. Moreover, since $A_p$ is a
linking matrix, each $A_{ii}$ is symmetric. By Lemma 2.7,
$\mbox{\it null}_p A_{ii}\ne 1$, hence, $\mbox{\it null}_p A_p\ne
1$. It follows from Lemma 2.9 that $H_1(M; {\bf Z}_p)\ne {\bf
Z}_p$. \hfill $\Box$

\medskip
To finish our proof of the main theorem of Section 2 we need the following
result (Corollary 2.3 in [Mu-2]).

\medskip
\noindent {\sc Proposition 2.11} (H. Murakami) {\it Fix an odd
prime $r$. For every connected, closed, oriented 3-manifold $M$,
there exist lens spaces $L(n_1, 1),  \dots , L(n_k, 1)$ with $n_i$
coprime to $r$ such that the connected sum $M \# L(n_1,1) \# \dots
\# L(n_k, 1)$ can be obtained by Dehn surgery on an algebraically
split link with integer framing.} \hfill $\Box$

\medskip

Now we are ready to prove Theorem 2.1 in full generality.

\medskip
\noindent
{\sc Proof of Theorem 2.1.}
Let $M$ be a $p$-periodic closed oriented 3-manifold and
$M_*=M/{\bf Z}_p$.
 By Proposition 2.11, there are integers $n_1,
\dots, n_k$ coprime to $p$ such that the connected sum
$\widetilde{ M_*} = M_*\# L(n_1,1) \# \dots \# L(n_k, 1)$ can be
obtained by Dehn surgery on an algebraically split framed link
$\widetilde{L_*}$.
Consider $\widetilde{ M} = M\# p L(n_1,1) \# \dots \# p L(n_k,
1)$.  Obviously, $\widetilde{ M}$ is $p$-periodic such that
$\widetilde{ M_*}= \widetilde{ M}/ {\bf Z}_p$. Moreover, it easily
follows from the proof of Theorem 1.1 that $\widetilde{ M}$ can be
obtained using Dehn surgery on an orbitally separated framed link
(see Remark 1.7).
Therefore, by Proposition 2.10,
$H_1(\widetilde{M}; {\bf Z}_p) \ne {\bf Z}_p$.
Since the numbers $n_1, \dots , n_k$ are coprime to $p$, we have
$H_1(L(n_i, 1), {\bf Z}_p)=0$ for every $i$. This implies that
$H_1(M; {\bf Z}_p)=H_1(\widetilde{M}; {\bf Z}_p) \ne {\bf Z}_p$.
 \hfill $\Box$

\bigskip
\noindent
{\bf 2.4. Orientation preserving ${\bf Z}_2$ actions.  Proof of Theorem 2.2.}

\medskip
\noindent
{\sc Remark 2.12.} Theorem 2.1 is not true in the case $p=2$. A simple
counterexample is $S^2\times S^1$.  It is interesting to
notice that $S^2\times S^1$ admits two different orientation
preserving actions of ${\bf Z}_2$ with the fixed point set being a
circle.  Indeed, let $H(1,1)$ be the negative Hopf link with
framing 1 on each component and let $H(-1, -1)$ be the negative
Hopf link with framing $-1$ on each component (see Fig. 3).  Dehn surgery on
each of these framed links produces $S^2\times S^1$.  Moreover,
permutation of the components defines ${\bf Z}_2$ actions on
$S^2\times S^1$ with a circle as the fixed point set.  These two
actions are different.  In the first case the orbit space of the
action $M_*$ is $S^2\times S^1$, and in the second case
$M_*=RP^3$.  Furthermore, one can see that in the first case ${\bf
Z}_2$ acts on ${\bf Z}=H_1(S^2\times S^1)$ trivially, and in the
second case it sends 1 to $-1$.  Both of the described actions
were studied in [P-3].

\ \\
\centerline{\psfig{figure=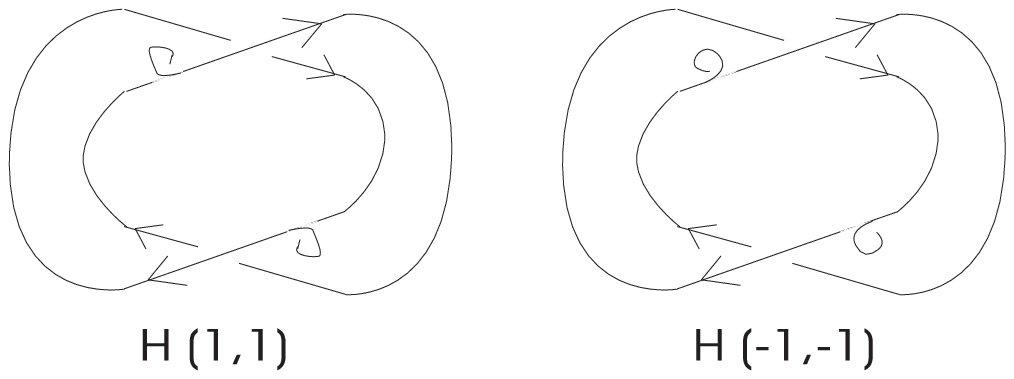}}
\begin{center}
Fig. 3
\end{center}

An interesting criterion for 2-periodic rational homology
3-spheres is provided by Theorem 2.2. Before we prove it, let us
recall that every finite abelian group can be uniquely decomposed
into a direct sum of cyclic groups whose orders are powers of
prime numbers.  Such decomposition will be called {\it canonical}.

\medskip
\noindent
{\sc Proof of Theorem 2.2.}
Let $M$ be a rational homology 3-sphere such that $H_1(M, {\bf
Z})$ does not have elements of order 16.  Assume that $M$ admits
an orientation-preserving action of ${\bf Z}_2$ such that the
fixed point set is a circle.  As before, let $M_*= M/{\bf Z}_2$
be the orbit space of the action.  By Theorem 1.1, $M$ can be
obtained using Dehn surgery on a strongly 2-periodic framed link
$L^2$ with the underlying link $L_*$.  By Lemma 1.4, $M_*$ is also
a rational homology 3-sphere.  Therefore, $M_*$ can be obtained by
surgery on an algebraically split framed link (see [Mu-1], [O-1]).
Thus, we may assume that $L_*$ is algebraically split (see Remark 1.7).
 By definition $L^2$ is orbitally
separated, and by Proposition 2.3 and Corollary 2.5 the linking
matrix of $L^2$ is block diagonal with every block being a
$2\times 2$ symmetric circulant matrix.  Recall that the linking
matrix of $L^2$ can be considered as a presentation matrix for the
abelian group  $H_1(M, {\bf Z})$.  Therefore, it is enough to show
that every finite abelian group $G$ presented by a matrix
 $A=\pmatrix{ a & b\cr
              b & a\cr}$ with $a, b\in {\bf Z}$ either has an
element of order 16 or there are even number of terms ${\bf
Z}_2$ and even number of terms ${\bf Z}_4$ in the canonical
decomposition of $G$, moreover it is possible that the canonical
decomposition of $G$ contains only one term of the form ${\bf
Z}_{2^t}$ for any $t\ge 3$. Denote by $g$ the greatest common
divisor of $a$ and $b$.  Since $G$ is finite, we have $det A\ne 0$
and $g>0$.  It is easy to see that $G$ can be presented by the
diagonal matrix  $\widetilde A=\pmatrix{ g & 0\cr
              0 & {{|det A|}\over g}\cr}$.
Therefore $G\simeq {\bf Z}_g \oplus {\bf Z}_{|det A|/g}$ (here by
${\bf Z}_1$ we mean the trivial group).

If $a$ and $b$ are both odd then $g$ is odd and  ${{|det A|}\over
g}$ is divisible by 8.

If $a$ and $b$ are both even then $g$ is even.  Let $t$ be the
power of 2 in the prime decomposition of $g$.  We have two
different cases: 1)  ${{|det A|}\over {g^2}}$  is odd. Then
$G\simeq {\bf Z}_{g/2^t}\oplus {\bf Z}_{2^t}\oplus {\bf Z}_{2^t} \oplus {\bf
Z}_{2l+1}$, where $2l+1={{|det A|}\over  {g 2^t}}$.
2) ${{|det A|}\over {g^2}}$ is even.  This means that both ${
a\over g}$ and $b\over g$ are odd and therefore  ${{|det A|}\over {g^2}}$
is divisible by $2^s, s\ge 3$, which implies that
$G\simeq {\bf Z}_{g/2^t}\oplus {\bf Z}_{2^t}\oplus {\bf Z}_{2^{s+t}} \oplus {\bf
Z}_{2h+1}$, where $2h+1={{|det A|}\over { g 2^{s+t}}}$.  In this case we have
an element of order $2^{s+t}$, which is at least 16.

We are left with the case when one of the numbers $a$ and $b$ is even
and the other one is odd.  But in this case $g$ is odd and $|det
A|$ is odd.

Therefore, if $G$ does not have elements of order 16 then the
number of terms ${\bf Z}_2$ and  the
number of terms ${\bf Z}_4$ in the canonical decomposition of $G$
are both even numbers.  To see that it is possible to have a
single ${\bf Z}_8$ term, consider the case $a=3$ and
$b=1$ (An example of a manifold with such first homology
is the lens space $L(8,3)$).
\hfill $\Box$

\medskip
\noindent
{\sc Corollary 2.13} {\it  Let $M$ be a rational homology 3-sphere such
that the group $H_1(M; {\bf Z})$ does not have elements of order
8.  If $M$ admits an orientation preserving action of ${\bf Z}_2$
with the fixed point set being a circle then $H_1(M; {\bf Z}_2)\simeq
{\bf Z}_2^{\oplus m}$ for some even integer $m$.}

przytyck@research.circ.gwu.edu

sokolov@gwu.edu

http://gwu.edu/\~{\hskip .1pt}sokolov/math\_page/mathematics.htm

\end{document}